\documentclass[12pt]{article}
\usepackage{graphicx}
\usepackage{amsmath}
\usepackage{amsfonts}
\usepackage{amssymb}
\newtheorem{theorem}{Theorem}

\newtheorem{corollary}[theorem]{Corollary}

\newtheorem{lemma}[theorem]{Lemma}

\newtheorem{proposition}[theorem]{Proposition}

\newenvironment{proof}[1][Proof]{\textbf{#1.} }{\ \rule{0.5em}{0.5em}}

\def\text{\hbox}

\def\a{\alpha}
\def\b{\beta}

\def\e{\varepsilon}
\def\f{\phi}

\def\p{\pi}

\def\s{\sigma}
\def\S{\Sigma}

\def\m{\mu}

\begin{document}

\title{Isomorphisms and Homeomorphisms of a Class of Graphs
and Spaces \footnote{Work performed under the auspices of
G.N.S.A.G.A. of C.N.R. of Italy; supported by the University of
Bologna, funds for selected research topics. The work of the first author is supported by Brazilina funds of
CNPq (Proc. 30.1103/80) and of Pronex (Project 107/97).}}
\author{S\'ostenes Lins \and Michele Mulazzani}

\maketitle

\begin{abstract}
{We solve the isomorphism problem for the whole class of
Lins-Mandel gems ({\it g\/}raphs {\it e\/}ncoded {\it m\/}anifold{\it s\/}).
We also present certain homeomorphisms of
branched cyclic coverings of two-bridge hyperbolic links. As a
consequence, we prove that, in in a wide subset of interesting
cases, the isomorphism conditions for Lins-Mandel gems are
equivalent to the homeomorphism conditions for the encoded
3-manifolds.
}
\end{abstract}

\section{Introduction}

Lins-Mandel spaces have been introduced in \cite{LM} as a direct
combinatorial generalization of lens spaces, by a four parameter
family of 4-coloured graphs. The encoded spaces $S(n,p,q,m)$ are
closed orientable 3-manifolds, possibly with isolated singular
points. This class of spaces have been extensively studied by
several researchers (see \cite{CG1}, \cite{CG2}, \cite{Ca1},
\cite{Ca2}, \cite{Ca3}, \cite{Ca4}, \cite{Ca5}, \cite{Do},
\cite{Gr}, \cite{JT}, \cite{LM}, \cite{Mu1}, \cite{Mu2}) and
appears to be fairly rich, since it contains several interesting
3-manifolds (see \cite{Li}), such as the Poincar\'e homology
sphere $S(5,3,2,1)\cong S(3,5,4,1)$ \cite{KS}, the Seifert-Weber
hyperbolic dodecahedron space $S(5,8,3,2)\cong S(5,8,3,3)$
\cite{ST}, the euclidean Hantzsche-Wendt manifold $S(3,5,2,1)$
\cite{Zi2}, the hyperbolic Fomenko-Matveev-Weeks manifold
$S(3,7,4,1)$ \cite{HW}, which is the hyperbolic 3-manifold with
the smallest known volume, and also an infinite family of
Brieskorn manifolds $M(n,p,2)\cong S(n,p,1,n-1)$ \cite{Mi}. The
necessary and sufficient conditions on the four parameters of a
Lins-Mandel coloured graph to encode a 3-manifold (i.e. to be a
3-gem) have been obtained in \cite{Mu1}. Moreover, \cite{Mu2}
shows that every Lins-Mandel manifold is a cyclic covering of
${\bf S}^3$, branched over a two-bridge knot or link.

In this paper we find the necessary and sufficient conditions for
the isomorphism of Lins-Mandel gems. Actually, this problem
has been studied in \cite{Ca5} for a particular subfamily of
graphs, called crystallizations, but its results appear to be
incorrect (see Remark 2).

We also present certain homeomorphisms of branched cyclic
coverings of two-bridge hyperbolic links, obtained from results of
Zimmermann \cite{Zi1} and Sakuma \cite{Sa}. As a consequence, we
prove that, in several interesting cases, the isomorphism
conditions for Lins-Mandel gems are equivalent to the
homeomorphism conditions for the encoded 3-manifolds.

\section{Isomorphisms of Lins-Mandel gems}

Regarding the theory of PL-manifolds represented by edge-coloured graphs, we
refer to \cite{FGG} and \cite{Li}. We recall that, in this theory, a
$(d+1)$-coloured graph encodes a $d$-dimensional pseudo-manifold.
When the encoded space is a manifold, the graph is called a {\it
$d$-gem}. Notice that every manifold admits representation by gems.

The family of Lins-Mandel 4-coloured graphs $${\cal
G}=\{G(n,p,q,m)\mid n,p\in {\bf Z}^{+},\,q\in {\bf Z}_{ 2p},\,m\in
{\bf Z}_n\}\,\footnote{By definition,
the third and the fourth coordinate of
$G(n,p,q,m)$ will be always considered mod $2p$ and mod $n$ respectively.}$$ has been defined in \cite{LM} by the following rules:
the set of vertices of $G=G(n,p,q,m)$ is $$V(G)={\bf Z}_n\times
{\bf Z}_{2p}$$ and the coloured edges are obtained by the
following four fixed-point-free involutions $\e_0,\e_1,\e_2,\e_3$
on $V(G)$

\vbox{ $$\e_0(i,j)=(i+m\m (j-q),1-j+2q),$$
$$\e_1(i,j)=(i,j-(-1)^j),$$
$$\e_2(i,j)=(i,j+(-1)^j),$$
$$\e_3(i,j)=(i+\m(j),1-j),$$ }

\noindent where $\m :{\bf Z}_{2p}\rightarrow\{-1,+1\}$ is the
function $$\m (j)=\left\{\begin{array}{lcc}
+1&\mbox{   if }&1\le j\le p\\-1&\mbox{otherwise}&\\
\end{array}\right..$$

For each $k\in \{0,1,2,3 \}$,  the vertices $ v,w\in V(G)$ are joined
by a $k$-edge if and only if $\e_{k}(v)=w$.

Roughly speaking, the graph $G(n,p,q,m)$ is constructed by taking
$n$ copies $C_i$ of a bicoloured cycle of length $2p$ involving
colours 1 and 2, for $i\in {\bf Z}_n$, so that $V(C_i)=\{(i,j)\mid
j\in {\bf Z}_{2p}\}$. The cycle $C_i$ is joined to the cycles
$C_{i\pm 1}$ by $3$-edges, and to the cycles $C_{i\pm m}$ by
$0$-edges.

Each $G(n,p,q,m)\in {\cal G}$ represents a 3-dimensional (possibly
singular) manifold $S(n,p,q,m)$ and $${\cal
S}=\{S(n,p,q,m)\mid n,p\in {\bf Z}^{+},\,q\in {\bf Z}_{ 2p},\,m\in
{\bf Z}_n\}$$
will be called the family of Lins-Mandel
spaces.
Since every $G(n,p,q,m)\in {\cal G}$ is bipartite, then
every $S(n,p,q,m)\in {\cal S}$ is an orientable (singular)
3-manifold. The spaces $S(n,p,q,m)$ and
$S(n,kp,kq,m)$ are homeomorphic \cite{Mu2}, therefore we shall
assume $\gcd(p,q)=1$ in the following, without loss of generality.

In most cases $S(n,p,q,m)$ is a genuine manifold (i.e.,
without singular points).

\begin{lemma} \label{Lemma 0}  \cite{Mu1} The graph $G(n,p,q,m)$
is a 3-gem and therefore $S(n,p,q,m)$ is a 3-manifold if and
only if either (i) $p$ is even or (ii) $p$ is odd and
$m=0,(-1)^q$. \end{lemma}


Figures 1 and 2 show two examples of graphs of the family,
both representing an interesting manifold. Observe that in the
figures the (missing) 0-edges connect vertices labelled with the
same letter.

\bigskip

\begin{figure}[ht]
 \begin{center}
 \includegraphics*[totalheight=8cm]{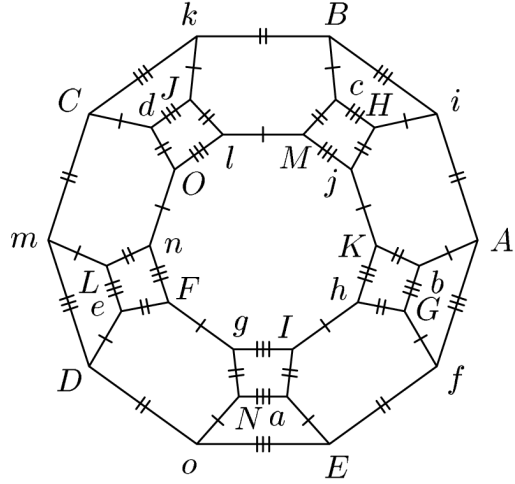}
 \end{center}
 \caption{$G(5,3,2,1)$, representing the Poincar\'e homology sphere.}

 \label{Fig. 1}

 \bigskip

\end{figure}






\bigskip

\begin{figure}[ht]
 \begin{center}
 \includegraphics*[totalheight=8cm]{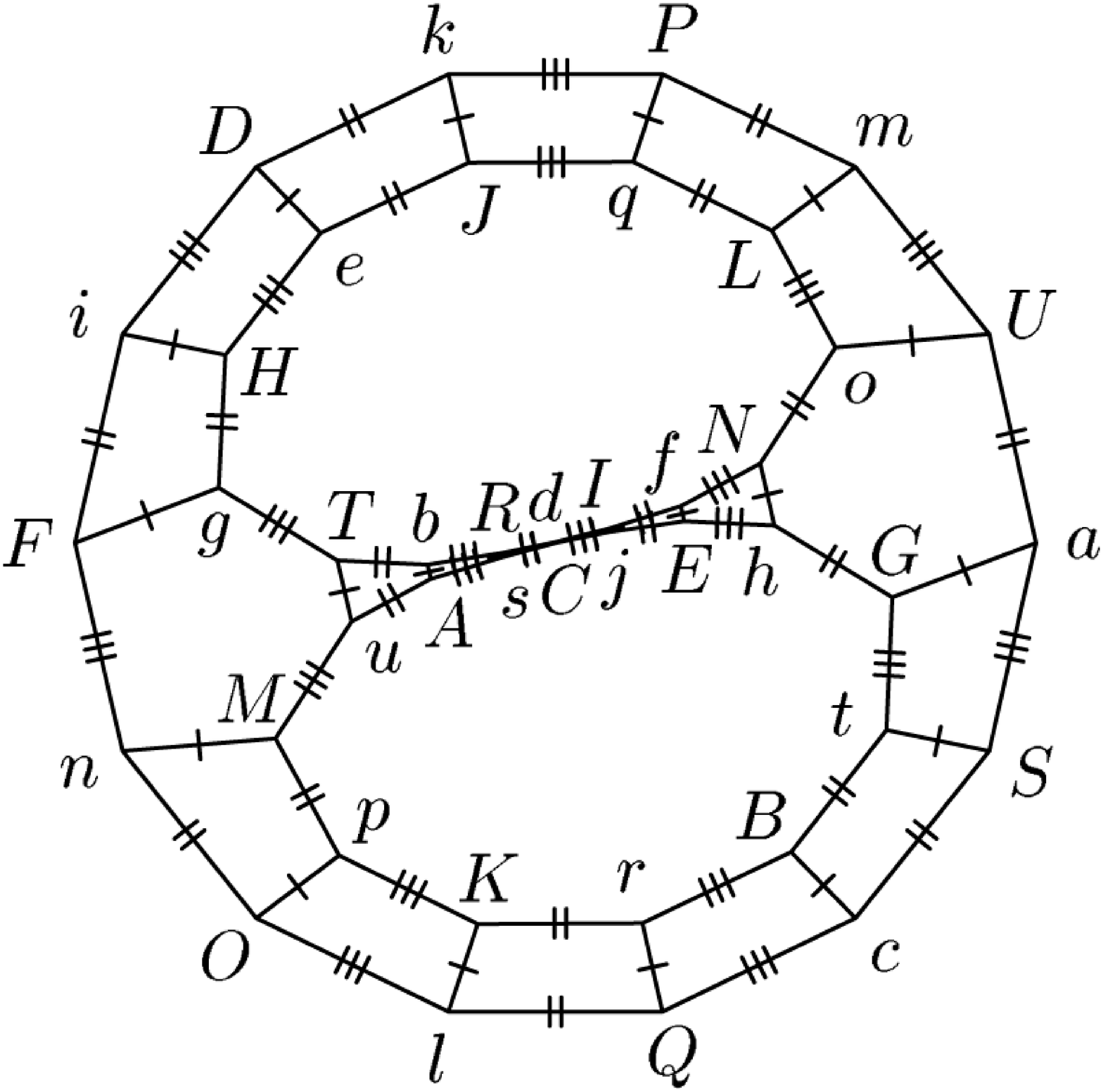}
 \end{center}
 \caption{$G(3,7,4,1)$, representing the Fomenko-Matveev-Weeks manifold.}

 \label{Fig. 3}

\bigskip

\end{figure}


\noindent {\bf Remark 1} The family of Lins-Mandel spaces has been
introduced as a combinatorial generalization of lens spaces. In
fact, for each $p,q>0$ such that $\gcd(p,q)=1$, the spaces
$S(2,p,q,1)$ and $S(p,2,1,q)$ are both homeomorphic to the lens
space $L(p,q)$. Moreover, $S(n,p,q,0)$ and $S(n,1,1,-1)$ are
homeomorphic to ${\bf S}^3$, for every $n,p>0$ and $q\in {\bf
Z}_{2p}$.

\medskip

Each graph $G(n,p,q,m)$ is defined by two different 4-tuples of
parameters, as stated by the following:

\begin{lemma} \label{Lemma 1} \cite{LM} The graphs $G(n,p,q,m)$ and
$G(n,p,q+p,-m)$ are equal. \end{lemma}

\begin{proof} Let $\e_0,\e_1,\e_2,\e_3$ (resp. $\e'_0,\e'_1,\e'_2,\e'_3$)
be the involutions defining $G(n,p,q,m)$ (resp. $G(n,p,q+p,-m)$).
Obviously, $\e'_1=\e_1$, $\e'_2=\e_2$ and $\e'_3=\e_3$. Moreover,
we have $\e_0'(i,j)=(i-m\m (j-q-p),1-j+2q-2p)=(i+m\m
(j-q),1-j+2q)= \e_0 (i,j)$. \footnote{Observe that $\m(j+p)=\m(j-p)=\m(1-j)=-\m(j)$,
for every $j\in{\bf Z}_{2p}$.}
\end{proof}

\medskip

Now we list the residues\footnote{If
$k',k''\in\{0,1,2,3\}$, a $\{k',k''\}$-residue is a bicoloured
cycle involving colours $k'$ and $k''$.} of $G(n,p,q,m)$ (see \cite{Mu1}).
When $p$ is even, the graph $G(n,p,q,m)$ contains:
\begin{itemize}
\item [-] $n$ $\{1,2\}$-residues of length $2p$,
\item [-] $n$ $\{0,3\}$-residues of length $2p$,
\item [-] $2$ $\{2,3\}$-residues of length $2n$ and $n(p-2)/2$ $\{2,3\}$-residues of length 4,
\item [-] $2\gcd(n,m)$ $\{0,1\}$-residues of length $2n/\gcd(n,m)$ and $n(p-2)/2$ $\{0,1\}$-residues of length 4,
\item [-] $np/2$ $\{1,3\}$-residues of length $4$,
\item [-] $np/2$ $\{0,2\}$-residues of length $4$.
\end{itemize}
On the other hand, when $p$ is odd, the graph $G(n,p,q,m)$ contains:
\begin{itemize}
\item [-] $n$ $\{1,2\}$-residues of length $2p$,
\item [-] $\gcd(n,m-(-1)^q)$ $\{0,3\}$-residues of length $2pn/\gcd(n,m-(-1)^q)$,
\item [-] $1$ $\{2,3\}$-residue of length $2n$ and $n(p-1)/2$ $\{2,3\}$-residues of length 4,
\item [-] $\gcd(n,m)$ $\{0,1\}$-residue of length $2n/\gcd(n,m)$ and $n(p-1)/2$ $\{0,1\}$-residues of length 4,
\item [-] $1$ $\{1,3\}$-residue of length $2n$ and $n(p-1)/2$ $\{1,3\}$-residues of length $4$,
\item [-] $\gcd(n,m)$ $\{0,2\}$-residue of length $2n/\gcd(n,m)$ and $n(p-1)/2$ $\{0,2\}$-residues of length $4$.
\end{itemize}

An isomorphism between two Lins-Mandel graphs $G$ and $G'$ is
uniquely defined by a pair $(f,\f)$, where $f:V(G)\to V(G')$ is a
bijection and $\f$ is a permutation of the colour set
$X=\{0,1,2,3\}$, such that $f\e_k=\e'_{\f (k)}f$, for each
$k\in\{0,1,2,3\}$. Actually, since $G$ and $G'$ are connected,
only $\f$ and the image $f(v)$ of a chosen vertex $v$ of $G$ are
required. Of course, isomorphic graphs encode homeomorphic spaces.

The next three lemmas give the main isomorphisms of Lins-Mandel
graphs. The proofs of these lemmas are rather technical and have
been included in the Appendix.

\begin{lemma} \label{Lemma A1} \cite{LM} $G(n,p,q,m)\cong G(n,p,-q,m)$.
\end{lemma}

\begin{lemma} \label{Lemma A2}
\begin{itemize}
\item [a)] If $p$ is even, then $G(n,p,q,m)\cong G(n,p,q^{-1},m)$.
\item [b')] If $p$ and $q$ are odd, then $G(n,p,q,-1)\cong G(n,p,q^{-1},-1)$.
\item [b'')] If $p$ is odd and $q$ is even, then $G(n,p,q,1)\cong G(n,p,(q+p)^{-1}+p,1)$.
\end{itemize}
\end{lemma}

\begin{lemma} \label{Lemma A3} \cite{Ca5} If $(n,m)=1$, then
$G(n,p,q,m)\cong G(n,p,q,m^{-1})$. \end{lemma}

From Lemmas \ref{Lemma 1}, \ref{Lemma A1}, \ref{Lemma A2} and
\ref{Lemma A3} we get:

\begin{corollary} \label{Corollary if}
\begin{itemize}
\item [a')] If $p$ is even, $\gcd(n,m)\ne 1$ and
$$\mbox  { either }\quad\left\{\begin{array}{lc}
q'=\pm q^{\pm 1}\\m'=m \\ \end{array}\right.
\quad \mbox  { or }\quad\left\{\begin{array}{lc}
q'=\pm q^{\pm 1}+p\\m'=-m\\ \end{array}\right.,$$
then $G(n,p,q',m')$ is isomorphic to $G(n,p,q,m)$.
\item [a'')] If $p$ is even, $\gcd(n,m)=1$ and
$$\mbox  { either }\quad\left\{\begin{array}{lc}
q'=\pm q^{\pm 1}\\m'=m^{\pm 1}\\ \end{array}\right.
\quad \mbox  { or }\quad\left\{\begin{array}{lc}
q'=\pm q^{\pm 1}+p\\m'=-m^{\pm 1}\\ \end{array}\right.,$$
then $G(n,p,q',m')$ is isomorphic to $G(n,p,q,m)$.
\item [b)] If $p$ is odd and $q'\equiv\pm q^{\pm 1}$ mod $p$, then
$G(n,p,q',(-1)^{q'})$ is isomorphic to $G(n,p,q,(-1)^{q})$.
\end{itemize}
\end{corollary}

\begin{proof} Statements a') and a'') follow from Lemmas \ref{Lemma 1},
\ref{Lemma A1}, \ref{Lemma A2} and \ref{Lemma A3}. As regards
part b), from the same lemmas we get:
\begin{itemize}
\item [b')] if $q$ is odd and $q'=\pm {q}^{\pm 1},\pm {q}^{\pm 1}+p$,
then $G(n,p,q',(-1)^{q'})\cong G(n,p,q,(-1)^{q})$;
\item [b'')] if $q$ is even and
$q'=\pm {(q+p)}^{\pm 1},\pm {(q+p)}^{\pm 1}+p$, then
$G(n,p,q',(-1)^{q'})\cong G(n,p,q,(-1)^{q})$.
\end{itemize}
It is easy to check that b')$\,+\,$b'') is equivalent to b).
\end{proof}

\medskip

Since for either $n\le 2$ or $p\le 2$ or $m=0$, the corresponding
space is trivial (${\bf S}^3$ or a lens space), we are mainly
interested in the cases $n,p\ge 3$ and $m\ne 0$.

\begin{lemma} \label{Lemma 2} \cite{Ca5} Let $n,p\ge 3$. If
$G(n',p',q',m')$ is isomorphic to $G(n,p,q,m)$, then $n'=n$ and
$p'=p$.
\end{lemma}

\begin{proof} As explained above, the $\{1,2\}$-residues of the graph $G(n,p,q,m)$
are exactly $n$ and they are all of length $2p$. When $n,p\ge 3$, the
same property holds for the $\{1,2\}$-residues (and possibly for
the $\{0,3\}$-residues) of $G(n',p',q',m')$ if and only if $n'=n$
and $p'=p$, while it does not hold for the other types of residues.
This proves the statement.
\end{proof}

\medskip

The previous result does not hold when either $n\le 2$ or $p\le 2$.
For example, $G(2,p,q,1)$ and $G(p,2,1,q)$ are isomorphic.

Corollary \ref{Corollary if} can be reversed when $n,p\ge 3$; the next
theorem completely describes the isomorphisms of these
graphs.

\begin{theorem} \label{Theorem 1} Assume $n,p\ge 3$.
\begin{itemize}
\item [a')] If $p$ is even and $\gcd(n,m)\ne 1$,
then $G(n',p',q',m')\cong G(n,p,q,m)$ if and only if
$$n'=n,\quad p'=p$$
and
$$\mbox { either }\quad\left\{\begin{array}{lc} q'=\pm
q^{\pm 1}\\m'=m \\
\end{array}\right.\quad \mbox  { or }\quad
\left\{\begin{array}{lc} q'=\pm q^{\pm 1}+p\\m'=-m\\
\end{array}\right..$$
\item [a'')] If $p$ is even and $\gcd(n,m)=1$, then
$G(n',p',q',m')\cong G(n,p,q,m)$ if and only if
$$n'=n,\quad p'=p$$
and
$$\mbox  { either }\quad\left\{\begin{array}{lc} q'=\pm q^{\pm
1}\\m'=m^{\pm 1}\\
\end{array}\right.\quad \mbox  { or }\quad
\left\{\begin{array}{lc} q'=\pm q^{\pm 1}+p\\m'=-m^{\pm 1}\\
\end{array}\right..$$
\item [b)] If $p$ is odd, then $G(n',p',q',(-1)^{q'})\cong
G(n,p,q,(-1)^q)$ if and only if $$n'=n,\quad p'=p\quad \mbox{and}
\quad q'\equiv\pm q^{\pm 1} \mbox{ mod } p .$$
\end{itemize}
\end{theorem}

\begin{proof}
The ``if'' part follows from Corollary \ref{Corollary if}. The
``only if'' part will be proved in the Appendix.
\end{proof}

\medskip

Observe that, for the second and the third parameter, the isomorphism conditions of part b) of Theorem
\ref{Theorem 1} are the same as the homeomorphism conditions for
lens spaces. This is not true for part a), since, in this case, the
situation is complicated by the presence of the additional parameter $m$.

\medskip

\noindent{\bf Remark 2} Proposition 4.1 of \cite{Ca5} states that,
when $\gcd(n,m)=1$, the graphs $G(n,p,q,m)$ and $G(n,p,q',m')$ are
isomorphic if and only if $q'\equiv\pm q^{\pm 1}$ mod $p$ and
$m'=\pm m^{\pm 1}$. This result is incorrect when $p$ is
even, since, for example, $H_1\big(S(3,4,1,1)\big)\cong {\bf
Z}_2\oplus{\bf Z}_6$ and
$H_1\big(S(3,4,5,1)\big)=H_1\big(S(3,4,1,2)\big)\cong {\bf Z}_3$
(see Appendix of \cite{LM}).

\medskip

Cases where $p$ is even are particularly interesting because the graph
always represents a manifold without any restriction on $m$. From
Theorem \ref{Theorem 1} we get:

\begin{corollary} \label{Corollary 1} Let $n,p,q$ be fixed,
with $n,p\ge 3$ and $p$ even. Then $G(n,p,q,m')\cong G(n,p,q,m)$
if and only if
\begin{itemize}
\item [(i)] $m'=m$, $\,$ when $\gcd(n,m)\ne 1$ and $q^2\ne p\pm 1$;
\item [(ii)] $m'=\pm m$, $\,$ when $\gcd(n,m)\ne 1$ and $q^2=p\pm 1$;
\item [(iii)] $m'=m^{\pm 1}$, $\,$ when $\gcd(n,m)=1$ and $q^2\ne p\pm 1$;
\item [(iv)] $m'=\pm m^{\pm 1}$, $\,$ when $\gcd(n,m)=1$ and $q^2=p\pm 1$.
\end{itemize}
\end{corollary}

\begin{proof} Since $p\ge 2$, $\gcd(p,q)=1$ and $q$ is odd, the condition $q=\pm q+p$ cannot
be satisfied and the condition $q=\pm q^{-1}+p$ is equivalent to
$q^2=\pm 1+p$. This proves the statement. \end{proof}


\section {Connections with branched cyclic coverings of two-bridge knots and links}

A $b$-fold branched cyclic covering of an oriented $\m$-component
link $L=\bigcup_{i=1}^{\m}L_i\subset{\bf S}^3$ is completely
determined (up to equivalence) by assigning to each component $L_i$
an integer $k_i\in{\bf Z}_b-\{0\}$, such that the set
$\{k_1,\ldots,k_{\m}\}$ generates the group ${\bf Z}_b$. The
monodromy associated to the covering sends each meridian of $L_i$
to the permutation $(1\,2\,\cdots\,b)^{k_i}\in\S_b$. By multiplying
each $k_i$ by the same invertible element $k$ of ${\bf Z}_b$, we
get an equivalent covering.

Following \cite{MM} we shall call a branched cyclic covering:
\begin{itemize}
\item [a)] {\it strictly-cyclic\/} if $k_i=k_j$, for every
$i,j\in\{1,\ldots,\m\}$,
\item [b)] {\it almost-strictly-cyclic\/} if $k_i=\pm k_j$, for every
$i,j\in\{1,\ldots,\m\}$,
\item [c)] {\it meridian-cyclic\/} if $\gcd(b,k_i)=1$, for every
$i\in\{1,\ldots,\m\}$,
\item [d)] {\it singly-cyclic\/} if $\gcd(b,k_i)=1$, for some
$i\in\{1,\ldots,\m\}$,
\item [e)] {\it monodromy-cyclic\/} if it is cyclic.
\end{itemize}

The following implications are straightforward: $$\text{ a)
}\Rightarrow\text{ b) }\Rightarrow\text{ c) }\Rightarrow\text{ d)
} \Rightarrow\text{ e) }.$$ Moreover, the five definitions are
equivalent when $L$ is a knot.

If $L$ has two components and the branched covering is
singly-cyclic, we can always suppose that $k_1=1$, up to equivalence
and renumbering of the components of $L$. Therefore, the covering
is completely determined by an integer $k=k_2\in{\bf Z}_b-\{0\}$.

Branched cyclic coverings of two-bridge knots and links are of
great interest, since a double branched covering of a two-bridge
knot or link is homeomorphic to a lens space (for notations and
properties of two-bridge knots and links we refer to \cite{BZ} and
\cite{Sc}). Let us denote by ${\bf b}(\a,\b)$ the two-bridge knot
or link of type $(\a,\b)$. It is well known that ${\bf b}(\a,\b)$
is a knot when $\a$ is odd and a two-component link when $\a$ is
even. Moreover, ${\bf b}(\a,\b)$ is hyperbolic if and only if it
is not toroidal (that is, $\b\not\equiv\pm 1$ mod $\a$). We denote by
$M_{b,k}(\a,\b)$ the $b$-fold singly-cyclic branched covering of
the link ${\bf b}(\a,\b)$ defined by $k$. Observe that the
branched covering $M_{b,k}(\a,\b)$ is strictly-cyclic if $k=1$,
almost-strictly-cyclic if $k=\pm 1$ and meridian-cyclic if
$\gcd(b,k)=1$.

From results of Zimmermann \cite{Zi1} and Sakuma \cite{Sa} it is
possible to obtain the homeomorphism conditions for these
manifolds, when the covering is meridian-cyclic and the branching
set is a hyperbolic link.

\begin{theorem} \label{Theorem 2}
Let $b,\a,\b$ be fixed, with $\a$ even and $\b\not\equiv\pm 1$ mod $\a$.
For $\gcd(b,k)=1$, the manifolds $M_{b,k'}(\a,\b)$ and
$M_{b,k}(\a,\b)$ are homeomorphic if and only if
\begin{itemize}
\item [(i)] $k'=k^{\pm 1}$, $\,$ when $\b^2\ne\a\pm 1$;
\item [(ii)] $k'=\pm k^{\pm 1}$, $\,$ when $\b^2=\a\pm1$.
\end{itemize}
\end{theorem}

\begin{proof}
Apply Theorem 1 of \cite{Zi1} (including the note
(a) of page 293) and Theorem 4.1 of \cite{Sa} (see tables of page 184).
\end{proof}

\medskip

Notice that, for the particular case of the Whitehead link ${\bf
b}(8,3)$, the previous result is contained in \cite{De} and
\cite{Zi1}.

\medskip

The graph $G(n,p,q,m)$ has a rotational cyclic symmetry of order
$n$, which sends each cycle $C_i$ onto $C_{i+1}$ (see details in
Lemma \ref{Lemma A1'}). As a direct consequence, the space
$S(n,p,q,m)$ also admits a cyclic symmetry. The next lemma states this
important property.

\begin{lemma} \label{Lemma Mich} \cite{Mu2}
\begin{itemize}
\item [a)] If $p$ is even and $m\ne 0$, then $S(n,p,q,m)$ is
homeomorphic to the singly-cyclic branched covering
$M_{n,-m}(p,q)$ of the two-bridge link ${\bf b}(p,q)$.
\item [b)] If $p$ is odd, then $S(n,p,q,(-1)^q)$ is the $n$-fold
branched cyclic covering of the two-bridge knot ${\bf b}(p,q)$.
\end{itemize}
\end{lemma}

As a consequence of the previous result, the geometric structure of
$S(n,p,q,m)$, when the branching set ${\bf b}(p,q)$ is a
hyperbolic knot or link, can be obtained from Thurston \cite{Th} and
Dunbar \cite{Du} results. Moreover, when the branching set is
toroidal and $m=(-1)^q$, then $S(n,p,q,m)$ turns out to be the
Brieskorn manifold $M(n,p,2)$ (see \cite{CG2} and \cite{Do}).
Thus, we have the following result for the geometric structure of
the manifold $S(n,p,q,m)$:

\begin{proposition} \label{Proposition 3} Let $S(n,p,q,m)$ be a
manifold.
\begin{itemize}
\item [1)] If $\gcd(n,m)=1$ and $q\not\equiv\pm 1$ mod $p$, then
$S(n,p,q,m)$ is hyperbolic for (i) $p=5$, $n\ge 4$ and (ii) $p\ne
5$, $n\ge 3$. Moreover, $S(3,5,3,-1)\cong S(3,5,2,1)$ is euclidean
and $S(2,p,q,1)$ is spherical for all $p,q$.
\item [2)] If $q\equiv\pm 1$ mod $p$, then $S(n,p,q,(-1)^q)$
is spherical for $n^{-1}+p^{-1}>1/2$, a Nil-manifold for
$n^{-1}+p^{-1}=1/2$ and a $\widetilde {SL}(2,\bf{R})$-manifold for
$n^{-1}+p^{-1}<1/2$.
\end{itemize}
\end{proposition}

\begin{proof} 1) See Theorem 3.1 and Remark 3.3 of \cite{HLM}. 2) See
\cite{Mi}. \end{proof}

\section {Homeomorphisms of Lins-Mandel manifolds}
Since isomorphic graphs encode homeomorphic spaces, from Corollary
\ref{Corollary if} we get the following homeomorphisms of
Lins-Mandel manifolds:

\begin{proposition} \label{Proposition if}
\begin{itemize}
\item [a')] If $p$ is even, $\gcd(n,m)\ne 1$ and $$\mbox  { either }\quad\left\{\begin{array}{lc}
q'=\pm q^{\pm 1}\\m'=m \\
\end{array}\right.\quad \mbox  { or }\quad
\left\{\begin{array}{lc} q'=\pm q^{\pm 1}+p\\m'=-m
\\
\end{array}\right.,$$ then $S(n,p,q',m')$ is homeomorphic to $S(n,p,q,m)$.
\item [a'')] If $p$ is even, $\gcd(n,m)=1$ and $$\mbox  { either }\quad\left\{\begin{array}{lc}
q'=\pm q^{\pm 1}\\m'=m^{\pm 1}\\
\end{array}\right.\quad \mbox  { or }\quad
\left\{\begin{array}{lc} q'=\pm q^{\pm 1}+p\\m'=-m^{\pm 1}\\
\end{array}\right.,$$ then $S(n,p,q',m')$ is homeomorphic to $S(n,p,q,m)$.
\item [b)] If $p$ is odd and $q'\equiv\pm q^{\pm 1}$ mod $p$, then
$S(n,p,q',(-1)^{q'})$ is homeomorphic to $S(n,p,q,(-1)^{q})$.
\end{itemize}
\end{proposition}

This proposition cannot be reversed, even when $n,p\ge 3$. For
example, \cite{Li} shows that $S(5,3,2,1)$ and $S(3,5,4,1)$
are both homeomorphic to the Poincar\'e homology sphere.

Theorem \ref{Theorem 2} and Lemma \ref{Lemma Mich} give us the
possibility of partially translating the combinatorial results of
Corollary \ref{Corollary 1} regarding graphs to the topological results for
manifolds. In fact, we have:

\begin{theorem} \label{Theorem 3} Let $n,p,q$ be fixed,
with $n,p\ge 3$, $p$ even and $q\not\equiv\pm 1$ mod $p$. For $\gcd(n,m)=1$,
the manifolds $S(n,p,q,m')$ and $S(n,p,q,m)$ are homeomorphic
if and only if
\begin{itemize}
\item [(i)] $m'=m^{\pm 1}$, $\,$ when $q^2\ne p\pm 1$;
\item [(ii)] $m'=\pm m^{\pm 1}$, $\,$ when $q^2=p\pm1$.
\end{itemize}
\end{theorem}

Observe that Theorem \ref{Theorem 3} is only a partial analogue of
Corollary \ref{Corollary 1}. So it is natural to state the
following:

\medskip

\noindent {\bf Conjecture}  Let $n,p,q$ be fixed, with $n,p\ge 3$,
$p$ even and $q\not\equiv\pm 1$ mod $p$. For $(n,m)\ne 1$, the manifolds
$S(n,p,q,m')$ and $S(n,p,q,m)$ are homeomorphic if and only if
\begin{itemize}
\item [(i)] $m'=m$, $\,$ when $q^2\ne p\pm 1$;
\item [(ii)] $m'=\pm m$, $\,$ when $q^2=p\pm1$.
\end{itemize}


\section{Appendix}
In this section we give the proofs of Lemmas \ref{Lemma A1},
\ref{Lemma A2} and \ref{Lemma A3} and complete the proof of Theorem
\ref{Theorem 1}.

\medskip

\begin{proof} {\bf of Lemma \ref{Lemma A1}} \cite{LM}
Let $\e_0,\e_1,\e_2,\e_3$ (resp. $\e'_0,\e'_1,\e'_2,\e'_3$) be the
involutions defining $G=G(n,p,q,m)$ (resp. $G'=G(n,p,-q,m)$). Moreover, let $\f_1\in\S_X$
and $f_1:{\bf Z}_n\times {\bf Z}_{2p}\to {\bf
Z}_n\times {\bf Z}_{2p}$ be the maps
$$\f_1=1,\qquad f_1(i,j)=(-i,1-j).$$
Since $f_1^2=1$, the map $f_1$ is a bijection. The pair $(f_1,\f_1)$ is
an isomorphism between the graphs $G(n,p,q,m)$ and $G(n,p,-q,m)$.
In fact, we get:
\begin{itemize}
\item [-] $f_1\e_0(i,j)=(-i-m\m(j-q),j-2q)$;
\item [-] $f_1\e_1(i,j)=(-i,1-j+(-1)^j)$;
\item [-] $f_1\e_2(i,j)=(-i,1-j-(-1)^j)$;
\item [-] $f_1\e_3(i,j)=(-i-\m(j),j)$;
\item [-] $\e'_0f_1(i,j)=(-i+m\m(1-j+q),j-2q)=(-i-m\m(j-q),j-2q)$;
\item [-] $\e'_1f_1(i,j)=(-i,1-j-(-1)^{1-j})=(-i,1-j+(-1)^j)$;
\item [-] $\e'_2f_1(i,j)=(-i,1-j+(-1)^{1-j})=(-i,1-j-(-1)^j)$;
\item [-] $\e'_3f_1(i,j)=(-i+\m(1-j),j)=(-i-\m(j),j)$.
\end{itemize}
\end{proof}

\medskip

In order to prove Lemma \ref{Lemma A2} we need some technical results.

\begin{lemma} \label{Lemma A0} \cite{Mu1} Let $\e_0,\e_1,\e_2,\e_3$ be the
involutions defining the graph $G(n,p,q,m)$.
\begin{itemize}
\item [a)] If $p$ is even, then: (i) $\e_3\e_1=\e_1\e_3$, (ii)
$\e_0\e_2=\e_2\e_0$, (iii) $\e_3\e_2(i,j)=\e_2\e_3(i,j)$ for every
$i$ and for every $j\ne 0,1,p,p+1$, (iv)
$\e_0\e_1(i,j)=\e_1\e_0(i,j)$ for every $i$ and for every $j\ne
q,q+1,q+p,q+p+1$.
\item [b')] If $p$ and $q$ are odd, then: (i) $\e_3\e_1(i,j)=\e_1\e_3(i,j)$
for every $i$ and for every $j\ne p,p+1$, (ii)
$\e_0\e_2(i,j)=\e_2\e_0(i,j)$ for every $i$ and for every $j\ne
q+p,q+p+1$, (iii) $\e_3\e_2(i,j)=\e_2\e_3(i,j)$ for every $i$ and
for every $j\ne 0,1$, (iv) $\e_0\e_1(i,j)=\e_1\e_0(i,j)$ for every
$i$ and for every $j\ne q,q+1$.
\item [b'')] If $p$ is odd and $q$ is even, then: (i)
$\e_3\e_1(i,j)=\e_1\e_3(i,j)$ for every $i$ and for every $j\ne
p,p+1$, (ii) $\e_0\e_2(i,j)=\e_2\e_0(i,j)$ for every $i$ and for
every $j\ne q,q+1$, (iii) $\e_3\e_2(i,j)=\e_2\e_3(i,j)$ for every
$i$ and for every $j\ne 0,1$, (iv) $\e_0\e_1(i,j)=\e_1\e_0(i,j)$
for every $i$ and for every $j\ne q+p,q+p+1$.
\end{itemize}
\end{lemma}

\begin{lemma} \label{Lemma A'0} \cite{Mu1} If $G(n,p,q,m)$ has $n$ $\{0,3\}$-residues,
then the $2p$ vertices of each $\{0,3\}$-residue have distinct
second coordinates.
\end{lemma}

\begin{lemma} \label{Lemma A''0} Let $(f,\f)$ be an isomorphism between $G=G(n,p,q,m)$
and $G'=G(n',p',q',m')$ and let $\e_0,\e_1,\e_2,\e_3$ (resp.
$\e'_0,\e'_1,\e'_2,\e'_3$) be the involutions defining $G$ (resp.
defining $G'$). If $\e'_{\f(k)}f(i,j)=f\e_k(i,j)$ and
$(i',j')=\e_k(i,j)$, then $\e'_{\f(k)}f(i',j')=f\e_k(i',j')$.
\end{lemma}

\begin{proof} We have $f\e_k(i',j')=f\e_k\e_k(i,j)=f(i,j)=
\e'_{\f(k)}\e'_{\f(k)}f(i,j)=\e'_{\f(k)}f\e_k(i,j)=\e'_{\f(k)}f(i',j')$.
\end{proof}

\medskip

Define $\p'':{\bf Z}_n\times {\bf Z}_{2p}\to{\bf Z}_{2p}$ as the
projection $\p''(i,j)=j$.

\medskip

\begin{proof} {\bf of Lemma \ref{Lemma A2}}
a) Let $\e_0,\e_1,\e_2,\e_3$ (resp. $\e'_0,\e'_1,\e'_2,\e'_3$) be
the involutions defining $G=G(n,p,q,m)$ (resp.
$G'=G(n,p,q^{-1},m)$). Moreover, let $\f_2\in\S_X$ and $f_2:{\bf
Z}_n\times {\bf Z}_{2p}\to {\bf Z}_n\times {\bf Z}_{2p}$ be the
maps $$\f_2=(0\,1)(2\,3),\qquad
f_2(i,j)=\left\{\begin{array}{lccc}(\e'_0
\e'_3)^{j/2}(-i,0)&\mbox{if} &j&\mbox{is even}\\
\e'_3(\e'_0\e'_3)^{(j-1)/2}(-i,0)&\mbox{if}&j&\mbox{is odd}\\
\end{array}\right..$$
Since $f_2$ sends the $\{1,2\}$-residues of $G$ injectively onto
the $\{0,3\}$-residues of $G'$, then it is a bijection. We claim
that $(f_2,\f_2)$ is an isomorphism between $G(n,p,q,m)$ and
$G(n,p,q^{-1},m)$. From the definition of $f_2$, we immediately
get $\e'_3f_2=f_2\e_2$ and $\e'_0f_2=f_2\e_1$. A direct
computation shows that $f_2(i,j)=(-i+h'_j,jq^{-1})$ if $j$ is even
and $f_2(i,j)=(-i+h''_j,1-(j-1)q^{-1})$ if $j$ is odd, where
$h'_j$ and $h''_j$ only depend on $j$. Since $\p''f_2(i,0)=0$,
$\p''f_2(i,1)=1$, $\p''f_2(i,p)=p$, $\p''f_2(i,p+1)=p+1$,
$\p''f_2(i,q)=q^{-1}$, $\p''f_2(i,q+1)=q^{-1}+1$,
$\p''f_2(i,q+p)=q^{-1}+p$ and $\p''f_2(i,q+p+1)=q^{-1}+p+1$, for
every $i\in{\bf Z}_n$, from Lemma \ref{Lemma A'0} we get:
$\p''f_2(i,j)\in\{0,1,p,p+1\}$ if and only if $j\in\{0,1,p,p+1\}$
and $\p''f_2(i,j)\in\{q^{-1},q^{-1}+1,q^{-1}+p,q^{-1}+p+1\}$ if
and only if $j\in\{q,q+1,q+p,q+p+1\}$. Now we show, by induction,
that (i') $\e'_2f_2(i,j)=f_2\e_3(i,j)$, for every $(i,j)\in{\bf
Z}_n\times {\bf Z}_{2p}$. By Lemma \ref{Lemma A''0} we only need
to prove (i') for $(i,1),(i,2),\ldots,(i,p)$. First of all, we have
$f_2\e_3(i,1)=f_2(i+1,0)=(-i-1,0)$ and
$\e'_2f_2(i,1)=\e'_2(-i-1,1)=(-i-1,0)$. Therefore, (i') holds in
$(i,1)$, for all $i$. Let us suppose that (i') holds in
$(i,1),(i,2),\ldots,(i,k)$, with $1\le k\le p-1$. From Lemma
\ref{Lemma A0} we get, if $k$ is odd,
$f_2\e_3(i,k+1)=f_2\e_3\e_1(i,k)=f_2\e_1\e_3(i,k)=\e'_0f_2\e_3(i,k)=
\e'_0\e'_2f_2(i,k)=\e'_2\e'_0f_2(i,k)=\e'_2f_2\e_1(i,k)=\e'_2f_2(i,k+1)$
and, if $k$ is even,
$f_2\e_3(i,k+1)=f_2\e_3\e_2(i,k)=f_2\e_2\e_3(i,k)=\e'_3f_2\e_3(i,k)=
\e'_3\e'_2f_2(i,k)=\e'_2\e'_3f_2(i,k)=\e'_2f_2\e_2(i,k)=\e'_2f_2(i,k+1)$.
Therefore, (i') holds in $(i,k+1)$, for all $i$. Then we show, by
induction, that (i'') $\e'_1f_2(i,j)=f_2\e_0(i,j)$, for every
$(i,j)\in{\bf Z}_n\times {\bf Z}_{2p}$. By Lemma \ref{Lemma A''0}
we only need to prove (i'') for $(i,q+1),(i,q+2),\ldots,(i,q+p)$.
First of all, we have
$f_2\e_0(i,q+1)=f_2(i+m,q)=(-i-m+h''_q,1-(q-1)q^{-1})=(-i-m+h''_q,q^{-1})$
and $\e'_1f_2(i,q+1)=\e'_1f_2\e_1(i,q)=\e'_1\e'_0f_2(i,q)=
\e'_1\e'_0(-i+h''_q,q^{-1})=\e'_1(-i-m+h''_q,1+q^{-1})=(-i-m+h''_q,q^{-1})$.
Therefore, (i'') holds in $(i,q+1)$, for all $i$. Let us suppose
that (i'') holds in $(i,q+1),(i,q+2),\ldots,(i,k)$ with $q+1\le
k\le q+p-1$. From Lemma \ref{Lemma A0} we get, if $k$ is odd,
$f_2\e_0(i,k+1)=f_2\e_0\e_1(i,k)=f_2\e_1\e_0(i,k)=\e'_0f_2\e_0(i,k)=
\e'_0\e'_1f_2(i,k)=\e'_1\e'_0f_2(i,k)=\e'_1f_2\e_1(i,k)=\e'_1f_2(i,k+1)$
and, if $k$ is even,
$f_2\e_0(i,k+1)=f_2\e_0\e_2(i,k)=f_2\e_2\e_0(i,k)=\e'_3f_2\e_0(i,k)=
\e'_3\e'_1f_2(i,k)=\e'_1\e'_3f_2(i,k)=\e'_1f_2\e_2(i,k)=\e'_1f_2(i,k+1)$.
Therefore (i') holds in $(i,k+1)$, for all $i$.

b') Let $\e_0,\e_1,\e_2,\e_3$ (resp. $\e'_0,\e'_1,\e'_2,\e'_3$) be
the involutions defining $G=G(n,p,q,-1)$ (resp.
$G'=G(n,p,q^{-1},-1)$). Moreover, let $\f_2\in\S_X$ and $f_2:{\bf
Z}_n\times {\bf Z}_{2p}\to {\bf Z}_n\times {\bf Z}_{2p}$ be the
maps $$\f_2=(0\,2)(1\,3),\qquad
f_2(i,j)=\left\{\begin{array}{lccc}(\e'_3
\e'_0)^{j/2}(-i,q^{-1}+1)&\mbox{if} &j&\mbox{is even}\\
\e'_0(\e'_3\e'_0)^{(j-1)/2}(-i,q^{-1}+1)&\mbox{if}&j&\mbox{is
odd}\\
\end{array}\right..$$
Since $f_2$ sends the $\{1,2\}$-residues of $G$ injectively onto
the $\{0,3\}$-residues of $G'$, then it is a bijection. We claim
that $(f_2,\f_2)$ is an isomorphism between $G(n,p,q,-1)$ and
$G(n,p,q^{-1},-1)$. From the definition of $f_2$, we immediately
get $\e'_3f_2=f_2\e_1$ and $\e'_0f_2=f_2\e_2$. A direct
computation shows that $f_2(i,j)=(-i+l'_j,1-(j-1)q^{-1})$ if $j$
is even and $f_2(i,j)=(-i+l''_j,jq^{-1})$ if $j$ is odd, where
$l'_j$ and $l''_j$ only depend on $j$. Since
$\p''f_2(i,0)=q^{-1}+1$, $\p''f_2(i,1)=q^{-1}$, $\p''f_2(i,p)=p$,
$\p''f_2(i,p+1)=p+1$, $\p''f_2(i,q)=1$, $\p''f_2(i,q+1)=0$,
$\p''f_2(i,q+p)=q^{-1}+p$ and $\p''f_2(i,q+p+1)=q^{-1}+p+1$, for
every $i\in{\bf Z}_n$, from Lemma \ref{Lemma A'0} we get:
$\p''f_2(i,j)\in\{0,1\}$ if and only if $j\in\{q,q+1\}$,
$\p''f_2(i,j)\in\{p,p+1\}$ if and only if $j\in\{p,p+1\}$,
$\p''f_2(i,j)\in\{q^{-1},q^{-1}+1\}$ if and only if $j\in\{0,1\}$
and $\p''f_2(i,j)\in\{q^{-1}+p,q^{-1}+p+1\}$ if and only if
$j\in\{q+p,q+p+1\}$. Now we show, by induction, that (i')
$\e'_1f_2(i,j)=f_2\e_3(i,j)$, for every $(i,j)\in{\bf Z}_n\times
{\bf Z}_{2p}$. By Lemma \ref{Lemma A''0} we only need to prove
(i') for $(i,1),(i,2)\ldots,(i,p)$. First of all, we have
$f_2\e_3(i,1)=f_2(i+1,0)=(-i-1,q^{-1}+1)$ and
$\e'_1f_2(i,1)=\e'_1(-i-1,q^{-1})=(-i-1,q^{-1}+1)$. Therefore, (i')
holds in $(i,1)$, for all $i$. Let us suppose that (i') holds in
$(i,1),(i,2),\ldots,(i,k)$ with $1\le k\le p-1$. From Lemma
\ref{Lemma A0} we get, if $k$ is odd,
$f_2\e_3(i,k+1)=f_2\e_3\e_1(i,k)=f_2\e_1\e_3(i,k)=\e'_3f_2\e_3(i,k)=
\e'_3\e'_1f_2(i,k)=\e'_1\e'_3f_2(i,k)=\e'_1f_2\e_1(i,k)=\e'_1f_2(i,k+1)$
and, if $k$ is even,
$f_2\e_3(i,k+1)=f_2\e_3\e_2(i,k)=f_2\e_2\e_3(i,k)=\e'_0f_2\e_3(i,k)=
\e'_0\e'_1f_2(i,k)=\e'_1\e'_0f_2(i,k)=\e'_1f_2\e_2(i,k)=\e'_1f_2(i,k+1)$.
Therefore, (i') holds in $(i,k+1)$, for all $i$. Then we show, by
induction, that (i'') $\e'_2f_2(i,j)=f_2\e_0(i,j)$, for every
$(i,j)\in{\bf Z}_n\times {\bf Z}_{2p}$. By Lemma \ref{Lemma A''0}
we only need to prove (i'') for $(i,q+1),(i,q+2),\ldots,(i,q+p)$.
First of all, we have $f_2\e_0(i,q+1)=f_2(i-1,q)=(-i+1+l''_q,1)$
and $\e'_2f_2(i,q+1)=\e'2f_2\e_1(i,q)=\e'_2\e'_3f_2(i,q)=
\e'_2\e'_3(-i+l''_q,1)=\e'_2(-i+1+l''_q,0)=(-i+1+l''_q,1)$.
Therefore, (i'') holds in $(i,q+1)$, for all $i$. Let us suppose
that (i'') holds in $(i,q+1), (i,q+2),\ldots,(i,k)$ with $q+1\le
k\le q+p-1$. From Lemma \ref{Lemma A0} we get, if $k$ is odd,
$f_2\e_0(i,k+1)=f_2\e_0\e_1(i,k)=f_2\e_1\e_0(i,k)=\e'_3f_2\e_0(i,k)=
\e'_3\e'_2f_2(i,k)=\e'_2\e'_3f_2(i,k)=\e'_2f_2\e_1(i,k)=\e'_2f_2(i,k+1)$
and, if $k$ is even,
$f_2\e_0(i,k+1)=f_2\e_0\e_2(i,k)=f_2\e_2\e_0(i,k)=\e'_0f_2\e_0(i,k)=
\e'_0\e'_2f_2(i,k)=\e'_2\e'_0f_2(i,k)=\e'_2f_2\e_2(i,k)=\e'_2f_2(i,k+1)$.
Therefore, (ii') holds in $(i,k+1)$, for all $i$.

b'') Follows directly from point b') and Lemma \ref{Lemma 1}.
\end{proof}

\medskip

\begin{proof} {\bf of Lemma \ref{Lemma A3}} \cite{Ca5}
Let $\e_0,\e_1,\e_2,\e_3$ (resp. $\e'_0,\e'_1,\e'_2,\e'_3$) be the
involutions defining $G=G(n,p,q,m)$ (resp. $G'=G(n,p,q,m^{-1})$).
Moreover, let $\f_3\in\S_X$ and $f_3:{\bf Z}_n\times {\bf Z}_{2p}\to {\bf
Z}_n\times {\bf Z}_{2p}$ be the maps
$$\f_3=\left\{\begin{array}{lccc}
(0\,3)(1\,2)&\mbox{if}&q&\mbox{is odd}\\
(0\,3)&\mbox{if}&q&\mbox{is even}\\
\end{array}\right.,\qquad f_3(i,j)=(-m^{-1}i,1+q-j).$$
It is easy to check that the map $g:{\bf Z}_n\times {\bf
Z}_{2p}\to {\bf Z}_n\times {\bf Z}_{2p}$, defined by
$g(i,j)=(-mi,1+q-j)$, is the inverse map of $f_3$; therefore, $f_3$ is a
bijection. The pair $(f_3,\f_3)$ is an isomorphism between $G$ and
$G'$. In fact we get:
\begin{itemize}
\item [-] $f_3\e_0(i,j)=(-m^{-1}i-\m(j-q),j-q)$;
\item [-] $f_3\e_1(i,j)=(-m^{-1}i,1+q-j+(-1)^j)$;
\item [-] $f_3\e_2(i,j)=(-m^{-1}i,1+q-j-(-1)^j)$;
\item [-] $f_3\e_3(i,j)=(-m^{-1}i-m^{-1}\m(j),q+j)$;
\item [-] $\e'_0f_3(i,j)=(-m^{-1}i+m^{-1}\m(1-j),q+j)=(-m^{-1}i-m^{-1}\m(j),q+j)$;
\item [-] $\e'_1f_3(i,j)=(-m^{-1}i,1+q-j-(-1)^{1+q-j})=
\left\{\begin{array}{lccc}
(-m^{-1}i,1+q-j+(-1)^j)&\mbox{if}&q&\mbox{is even}\\
(-m^{-1}i,1+q-j-(-1)^j)&\mbox{if}&q&\mbox{is odd}\\
\end{array}\right.$;
\item [-] $\e'_2f_3(i,j)=(-m^{-1}i,1+q-j+(-1)^{1+q-j})=
\left\{\begin{array}{lccc}
(-m^{-1}i,1+q-j-(-1)^j)&\mbox{if}&q&\mbox{is even}\\
(-m^{-1}i,1+q-j+(-1)^j)&\mbox{if}&q&\mbox{is odd}\\
\end{array}\right.$;
\item [-] $\e'_3f_3(i,j)=(-m^{-1}i+\m(1+q-j),j-q)=(-m^{-1}i-\m(j-q),j-q)$.
\end{itemize}
\end{proof}

\medskip

In order to prove Theorem \ref{Theorem 1} we need some preparatory
results.

\begin{lemma} \label{Lemma A1'} Let $\s\in\S_X$ and $r,s:{\bf Z}_n\times {\bf Z}_{2p}\to {\bf
Z}_n\times {\bf Z}_{2p}$ be the maps

$$\s=\left\{\begin{array}{lccc}1&\mbox{if}
&p&\mbox{is even}\\(1\,2)&\mbox{if}&p&\mbox{is odd}\\
\end{array}\right.,\quad r(i,j)=(i+1,j),\quad s(i,j)=(-i,p+j).$$
Then $(r,1)$ and $(s,\s)$ are automorphisms of $G(n,p,q,m)$.
\end{lemma}
\begin{proof} Since $r^n=1=s^2$,
both $r$ and $s$ are bijections. Then we have:
\begin{itemize}
\item [-] $\e_0r(i,j)=(i+1+m\m(j-q),1-j+2q)=r\e_0(i,j)$,
\item [-] $\e_1r(i,j)=(i+1,j-(-1)^{j})=r\e_1(i,j)$,
\item [-] $\e_2r(i,j)=(i+1,j+(-1)^{j})=r\e_2(i,j)$,
\item [-] $\e_3r(i,j)=(i+1+\m(j),1-j)=r\e_3(i,j)$.
\item [-] $\e_0s(i,j)=(-i-m\m(j-q),1+p-j+2q)$,
\item [-] $\e_1s(i,j)=(-i,p+j-(-1)^{p+j})=
\left\{\begin{array}{lccc} (-i,p+j-(-1)^j)&\mbox{if}&p&\mbox{is
even}\\ (-i,p+j+(-1)^j)&\mbox{if}&p&\mbox{is odd}\\
\end{array}\right.$;
\item [-] $\e_2s(i,j)=(-i,p+j+(-1)^{p+j})=\left\{\begin{array}{lccc}
(-i,p+j+(-1)^j)&\mbox{if}&p&\mbox{is
even}\\ (-i,p+j-(-1)^j)&\mbox{if}&p&\mbox{is odd}\\
\end{array}\right.$;
\item [-] $\e_3s(i,j)=(-i-\m(j),1-j+p)$,
\item [-] $s\e_0(i,j)=(-i-m\m(j-q),1+p-j+2q)$,
\item [-] $s\e_1(i,j)=(-i,j+p-(-1)^{j}$,
\item [-] $s\e_2(i,j)=(-i,j+p+(-1)^{j}$,
\item [-] $s\e_3(i,j)=(-i-\m(j),1-j+p)$.
\end{itemize} \end{proof}

\begin{lemma} \label{Lemma A1''} Let $(f,\f)$ be
an isomorphism between the graphs $G(n,p,q,m)$ and
$G(n,p,q',m')$, with $n,p\ge 3$.
\begin{itemize}
\item [a')] If $p$ is even and $\gcd(n,m)\ne 1$ then $\f\in\{1,(0\,1)(2\,3)\}$
and $\p''f(0,0)\in\{0,1,p,p+1\}$.
\item [a'')] If $p$ is even and $\gcd(n,m)=1$ then
$\f\in\{1,(0\,1)(2\,3),(0\,3)(1\,2),(0\,2)(1\,3)\}$ and
(i) $\p''f(0,0)\in\{0,1,p,p+1\}$ if $\f\in\{1,(0\,1)(2\,3)\}$,
(ii) $\p''f(0,0)\in\{q',q'+1,q'+p,q'+p+1\}$ if
$\f\in\{(0\,3)(1\,2),(0\,2)(1\,3)\}$.
\item [b)] If $p$ is odd then
$\f\in\{1,(1\,2),(0\,3),(0\,3)(1\,2),(0\,1)(2\,3),(0\,2)(1\,3),
(0\,1\,3\,2),(0\,2\,3\,1)\}$ and
$\p''f(0,0)\in\{0,1,q',q'+1,p,p+1,q'+p,q'+p+1\}$.
\end{itemize}
\end{lemma}

\begin{proof} First of all, observe that the vertex $(0,0)$ of $G=G(n,p,q,m)$ lies in a
$\{2,3\}$-residue of length $2n$. All the $\{1,2\}$-residues
(resp. the $\{1,2\}$-residues) of $G$ are mapped by any
isomorphism either to the $\{1,2\}$-residues or to the
$\{0,3\}$-residues of $G'=G(n,p,q',m')$.

a') Each of the two $\{2,3\}$-residues of length $2n$ is mapped to
a $\{2,3\}$-residue (where each vertex has second coordinate
$j\in\{0,1,p,p+1\}$) and each of the two $\{0,1\}$-residues of
length $2n/\gcd(n,m)$ is mapped to a $\{0,1\}$-residue.

a'') Each of the two $\{2,3\}$-residues (resp. the two
$\{0,1\}$-residues) of length $2n$ is mapped either to a
$\{2,3\}$-residue (where each vertex has second coordinates
$j\in\{0,1,p,p+1\}$) or to a $\{0,1\}$-residue (where each vertex
has second coordinates $j\in\{q',q'+1,q'+p,q'+p+1\}$).

b) The $\{2,3\}$-residue, the $\{1,3\}$-residue, the
$\{0,1\}$-residue and the $\{0,2\}$-residue of length $2n$ are
mapped to either the $\{2,3\}$-residue or the $\{1,3\}$-residue or
the $\{0,1\}$-residue or the $\{0,2\}$-residue. Since each vertex
of these four residues has second coordinate
$j\in\{0,1,q',q'+1,p,p+1,q'+p,q'+p+1\}$), the statement holds.
\end{proof}

\medskip

\begin{proof} {\bf of Theorem \ref{Theorem 1}} {\it (``only if''
part:)} From Lemma \ref{Lemma 2} we get $n'=n$ and $p'=p$.
Let $\e_0,\e_1,\e_2,\e_3$
(resp. $\e'_0,\e'_1,\e'_2,\e'_3$) be the involutions defining
either $G=G(n,p,q,m)$ or $G=G(n,p,q,(-1)^{q})$ (resp. either
$G'=G(n,p,q',m')$ or $G'=G(n,p,q',(-1)^{q'})$) and let
$f_1,f_2,f_3$ be the bijections defined in Lemma \ref{Lemma A1},
\ref{Lemma A2} and \ref{Lemma A3} respectively. Moreover, observe
that the pairs $(r,1)$ and $(s,\s)$ defined in Lemma \ref{Lemma
A1'} are automorphisms of both $G$ and $G'$.

a') If $(f,\f)$ is an isomorphism between $G$ and $G'$ then,
by Lemma \ref{Lemma A1''},
$\f\in\{1,(0\,1)(2\,3)\}$ and, up to the action of $r$ and $s$, we
can suppose $f(0,0)=(0,j)$ with $j\in\{0,1\}$. Then we have the
following four cases:
\begin{itemize}
\item [(i)] If $f(0,0)=(0,0)$ and $\f=1$, then $f=1$ and
$(f,\f)$ is an automorphism of $G(n,p,q,m)$.
\item [(ii)] If $f(0,0)=(0,0)$ and $\f=(0\,1)(2\,3)$, then $f=f_2$ and
$(f,\f)$ is an isomorphism between $G(n,p,q,m)$ and $G(n,p,q^{-1},m)$.
\item [(iii)] If $f(0,0)=(0,1)$ and $\f=1$, then $f=f_1$ and
$(f,\f)$ is an isomorphism between $G(n,p,q,m)$ and $G(n,p,-q,m)$.
\item [(iv)] If $f(0,0)=(0,1)$ and $\f=(0\,1)(2\,3)$, then $f=f_1f_2$ and
$(f,\f)$ is an isomorphism between $G(n,p,q,m)$ and $G(n,p,-q^{-1},m)$.
\end{itemize}

a'') If $(f,\f)$ is an isomorphism between $G$ and $G'$ then, by Lemma
\ref{Lemma A1''}, $\f\in\{1,(0\,3)(1\,2),(0\,1)(2\,3),(0\,2)(1\,3)\}$ and,
up to the action of $r$ and $s$, we can suppose $f(0,0)=(0,j)$, with
$j\in\{0,1\}$ if $\f\in\{1,(0\,1)(2\,3)\}$ and with
$j\in\{q',q'+1\}$ if $\f\in\{(0\,3)(1\,2),(0\,2)(1\,3)\}$.
Then we have the four cases of the previous point and the following four cases:
(i') $f(0,0)=(0,1+q')$ and $\f=(0\,3)(1\,2)$;
(ii') $f(0,0)=(0,1+q')$ and $\f=(0\,2)(1\,3)$;
(iii') $f(0,0)=(0,q')$ and $\f=(0\,3)(1\,2)$ and
(iv') $f(0,0)=(0,q')$ and $\f=(0\,2)(1\,3)$.
In (i'), if $j$ is even (resp. if $j$ is odd) we get
$f(i,j)=f(\e_1\e_2)^{j/2}(\e_3\e_2)^i(0,0)=(\e'_2\e'_1)^{j/2}(\e'_0\e'_1)^if(0,0)=(\e'_2\e'_1)^{j/2}(-m'i,1+q')=
(-m'i,1+q'-j)$ (resp.
$f(i,j)=f\e_2(\e_1\e_2)^{(j-1)/2}(\e_3\e_2)^i(0,0)=
\e'_1(\e'_2\e'_1)^{(j-1)/2}(\e'_0\e'_1)^if(0,0)=\e'_1(\e'_2\e'_1)^{(j-1)/2}(-m'i,1+q')=(-m'i,1+q'-j)$).
Since $f\e_0(i,j)=(-m'i-m'm\m(j-q),j+q'-2q)$ and
$\e'_3f(i,j)=(-m'i-\m(j-q'),j-q')$, we get $q'=q$, $m'=m^{-1}$,
and therefore $f=f_3$. Case (ii') is a composition of cases (i')
and (ii), case (iii') is a composition of cases (i') and (iii),
case (iv') is a composition of cases (i') and (iv). Therefore, we
have:
\begin{itemize}
\item [(i')] If $f(0,0)=(0,1+q)$ and $\f=(0\,3)(1\,2)$, then $f=f_3$ and
$(f,\f)$ is an isomorphism between $G(n,p,q,m)$ and $G(n,p,q,m^{-1})$.
\item [(ii')] If $f(0,0)=(0,1+q)$ and $\f=(0\,2)(1\,3)$, then $f=f_3f_2$ and
$(f,\f)$ is an isomorphism between $G(n,p,q,m)$ and $G(n,p,q^{-1},m^{-1})$.
\item [(iii')] If $f(0,0)=(0,q)$ and $\f=(0\,3)(1\,2)$, then $f$ is the map $f_3f_1$ and
$(f,\f)$ is an isomorphism between $G(n,p,q,m)$ and $G(n,p,-q,m^{-1})$.
\item [(iv')] If $f(0,0)=(0,q)$ and $\f=(0\,2)(1\,3)$, then $f$ is the map $f_3f_1f_2$ and
$(f,\f)$ is an isomorphism between $G(n,p,q,m)$ and $G(n,p,-q^{-1},m^{-1})$.
\end{itemize}

b) The pair $(f_3,(1\,2))$ is an automorphism of $G'$. Therefore,
if $(f,\f)$ is an isomorphism between $G$ and $G'$ then, up to the
action of $r$, $s$ and $f_3$, we can suppose that $f(0,0)=(0,j)$, with
$j\in\{0,1\}$. By this assumption, $f$ sends the $\{2,3\}$-residue
of $G$ onto the $\{2,3\}$-residue of $G'$ and therefore either
$\f=1$ or $\f=(0\,1)(2\,3)$. We have four cases:
\begin{itemize}
\item [(i'')] If $f(0,0)=(0,0)$ and $\f=1$, then $f=1$ and
$(f,\f)$ is an automorphism of $G(n,p,q,(-1)^q)$.
\item [(ii'')] If $f(0,0)=(0,0)$ and $\f=(0\,1)(2\,3)$, then $f=f_2$ and
$(f,\f)$ is an isomorphism between $G(n,p,q,-1)$ and
$G(n,p,q^{-1},-1)$ (resp. between $G(n,p,q,1)$ and
$G(n,p,(q+p)^{-1}+p,1)$) if $q$ is odd (resp. even).
\item [(iii'')] If $f(0,0)=(0,1)$ and $\f=1$, then $f=f_1$ and
$(f,\f)$ is an isomorphism between $G(n,p,q,(-1)^q)$ and $G(n,p,-q,(-1)^q)$.
\item [(iv'')] If $f(0,0)=(0,1)$ and $\f=(0\,1)(2\,3)$, then $f=f_1f_2$ and
$(f,\f)$ is an isomorphism between $G(n,p,q,-1)$ and
$G(n,p,-q^{-1},-1)$ (resp.  between $G(n,p,q,1)$ and
$G(n,p,-(q+p)^{-1}+p,1)$) if $q$ is odd (resp. even).
\end{itemize}
Therefore, if $G(n,p,q',(-1)^{q'})$ is isomorphic to
$G(n,p,q,(-1)^q)$, then ($\ast$) $q'=\pm q^{\pm 1}$ when $q$ is
odd and ($\ast\ast$) $q'=\pm(q+p)^{\pm 1}+p$ when $q$ is even.
Since ($\ast$) $+$ ($\ast\ast$) is equivalent to $q'\equiv\pm
q^{\pm 1}$ mod $p$, our proof is completed.
\end{proof}

\vspace{15 pt} {S\'OSTENES LINS, Departamento de Matem\`atica,
Universidade Federal de Pernambuco, Recife-PE, BRAZIL. E-mail:
sostenes@dmat.ufpe.br}

\vspace{15 pt} {MICHELE MULAZZANI, Dipartimento di Matematica and C.I.R.A.M.,
Universit\`a di Bologna, I-40127 Bologna, ITALY. E-mail: mulazza@dm.unibo.it}

\end{document}